\documentclass[letterpaper, 10 pt, conference]{ieeeconf}  

\IEEEoverridecommandlockouts                              

\makeatletter
\newcommand\fs@spaceruled{\def\@fs@cfont{\bfseries}\let\@fs@capt\floatc@ruled
	\def\@fs@pre{\vspace{0.4em}\hrule height.8pt depth0pt \kern2pt}%
	\def\@fs@post{\kern2pt\hrule\relax}%
	\def\@fs@mid{\kern2pt\hrule\kern2pt}%
	\let\@fs@iftopcapt\iftrue}
\makeatother

\usepackage{amsmath, amssymb, bbm, xspace}
\usepackage{amsfonts,mathrsfs}
\usepackage{mathtools}
\usepackage{color}
\usepackage{graphicx}
\usepackage{epsfig}
\usepackage{algorithm}
\usepackage{algorithmicx}
\usepackage[acronym]{glossaries}
\usepackage{subfigure}
\usepackage{cancel}
\usepackage{empheq}
\usepackage{placeins}
\let\labelindent\relax
\usepackage{enumitem}
\usepackage{textgreek}
\usepackage{url}
\def\be{\begin{equation}}
\def\ee{\end{equation}}
\def\ba{\begin{align}}
\def\ea{\end{align}}
\def\bas{\begin{align*}}
\def\eas{\end{align*}}
\newcommand{\bs}{\boldsymbol}
\newcommand{\mc}{\mathcal}

\newcommand{\blue}{\textcolor{black}}

\definecolor{darkgreen}{rgb}{0.0, 0.5, 0.0}

\renewcommand{\emph}{\textit}
\newcommand\numberthis{\addtocounter{equation}{1}\tag{\theequation}}

\makeglossaries
\newacronym{GNEP}{GNEP}{generalized Nash equilibrium problem}
\newacronym{NE}{NE}{Nash equilibrium}
\newacronym{NEP}{NEP}{Nash equilibrium problem}
\newacronym{GNE}{GNE}{generalized Nash equilibrium}
\newacronym{v-GNE}{v-GNE}{variational \gls{GNE}}
\newacronym{ISS}{ISS}{Input-to-state-stable}
\newacronym{PPA}{PPA}{proximal-point algorithm}
\newacronym{PPPA}{PPPA}{preconditioned \gls{PPA}}
\newacronym{VI}{VI}{variational inequality}
\newacronym{GAE}{GAE}{generalized aggregative equilibrium}
\newacronym{v-GAE}{v-GAE}{variational \gls{GAE}}
\newacronym{KKT}{KKT}{Karush--Kuhn--Tucker}
\newacronym{FQNE}{FQNE}{firmly quasinonexpansive}
\newacronym{FNE}{FNE}{firmly nonexpansive}
\newacronym{ADMM}{ADMM}{alternating direction method of multipliers}
\newacronym{MPMM}{MPMM}{modified proximal method of multipliers}
\newacronym{OPF}{OPF}{optimal power flow}
\newacronym{OPFP}{OPFP}{optimal power flow problem}
\newacronym{NUM}{NUM}{network utility maximization}

\newcommand{\0}{\bs 0}
\def\1{{\bs 1}}

\def\argmin{\mathop{\rm argmin}}

\newcommand{\col}{\mathrm{col}}

\def\diag{\mathop{\hbox{\rm diag}}}

\def\spose#1{\hbox to 0pt{#1\hss}}

\newcommand{\proj}{\mathrm{proj}}

\def\dom{\operatorname{dom}}

\def\R{\mathbb{R}}
\def\N{\mathbb{N}}

\def\I{\mc{I}}
\def\i{{i\in\mc{I}}}
\def\neigi{\mc{N}_i}

\usepackage{mathbbol}
\DeclareSymbolFontAlphabet{\mathbbm}{bbold}
\DeclareSymbolFontAlphabet{\mathbb}{AMSb}%
\def\h{\mathbbm{x}^\star}

\DeclareSymbolFont{myletters}{OML}{ztmcm}{m}{it}
\DeclareMathSymbol{\uplambda}{\mathord}{myletters}{"15}

\usepackage{xcolor}

\def\QEDhereeqn{\eqno\let\eqno\relax\let\leqno\relax\let\veqno\relax\hbox{\QED}}
\def\QEDopenhereeqn{\eqno\let\eqno\relax\let\leqno\relax\let\veqno\relax\hbox{\QEDopen}}

\renewcommand{\indent}{\noindent\hspace*{8pt}\ignorespaces}

\makeatletter \let\cl@part\relax \makeatother
\usepackage{nohyperref}
\usepackage[capitalize]{cleveref}
\def\k{{k \in \N}}  
\usepackage{letltxmacro}

\crefname{thm}{Theorem}{Theorem}
\crefname{lem}{Lemma}{Lemma}
\crefname{cor}{Corollary}{Corollary}
\crefname{claim}{Claim}{Claim}
\crefname{axiom}{Axiom}{Axiom}
\crefname{conj}{Conjecture}{Conjecture}
\crefname{fact}{Fact}{Fact}
\crefname{hypo}{Hypothesis}{Hypothesis}
\crefname{assum}{Assumption}{Assumption}
\crefname{prop}{Proposition}{Proposition}
\crefname{crit}{Criterion}{Criterion}
\crefname{standing}{Standing Assumption}{Standing Assumption}
\crefname{defn}{Definition}{Definition}
\crefname{exmp}{Example}{Example}
\crefname{rem}{Remark}{Remark}
\crefname{prob}{Problem}{Problem}
\crefname{prin}{Principle}{Principle}
\crefname{alg}{Algorithm}{Algorithm}
\crefname{figure}{Figure}{Figures}
\crefname{assumption}{Assumption}{Assumptions}

\crefname{thmlisti}{Theorem}{Theorems}
\crefname{lemlisti}{Lemma}{Lemma}
\crefname{asmlisti}{Assumption}{Assumption}

\newlist{thmlist}{enumerate}{1}
\setlist[thmlist]{label=(\roman{thmlisti}), ref=\thethm(\roman{thmlisti}),noitemsep}
\newlist{lemlist}{enumerate}{1}
\setlist[lemlist]{label=(\roman{lemlisti}), ref=\thelem(\roman{lemlisti}),noitemsep}
\newlist{asmlist}{enumerate}{1}
\setlist[asmlist]{label=(\roman{asmlisti}), ref=\theassumption(\roman{asmlisti}),noitemsep,nosep,leftmargin=*} 

\newtheorem{lemma}{Lemma}
\newtheorem{theorem}{Theorem}
\newtheorem{remark}{Remark}
\newtheorem{assumption}{Assumption}

\usepackage{stackengine}
\stackMath
\newcommand\tsup[2][2]{%
	\def\useanchorwidth{T}%
	\ifnum#1>1%
	\stackon[-.5pt]{\tsup[\numexpr#1-1\relax]{#2}}{\scriptscriptstyle\sim}%
	\else%
	\stackon[.5pt]{#2}{\scriptscriptstyle\sim}%
	\fi%
}

\title{\LARGE \bf
The distributed  dual ascent algorithm is robust to asynchrony
}

\author{Mattia Bianchi, Wicak Ananduta,  and Sergio Grammatico
\thanks{The authors are with the Delft Center for Systems and Control (DCSC), TU Delft, The Netherlands.
	E-mail addresses: \texttt{\{m.bianchi, w.ananduta, s.grammatico\}@tudelft.nl}. This work was partially supported by NWO under research project OMEGA (grant n. 613.001.702) and by the ERC under research project COSMOS (802348).} 
}

\begin{document}

\maketitle
\thispagestyle{empty}
\pagestyle{empty}

\begin{abstract}
The distributed dual ascent is an established algorithm to solve strongly convex multi-agent optimization problems with separable  cost functions, in the presence of coupling constraints. In this paper, we study its asynchronous counterpart.
Specifically, we assume that each agent only relies on the outdated information received from some  neighbors. Differently from the existing randomized and dual block-coordinate schemes, we show convergence under heterogeneous delays, communication and update frequencies. Consequently, our  asynchronous dual ascent algorithm can be implemented without requiring any coordination between the agents.
	
\end{abstract}

\section{Introduction}\label{sec:introduction}

Distributed multi-agent  optimization is 
well suited for
modern large-scale 
and data-intensive
problems,
where the volume and spatial scattering of the information render centralized processing and storage  inefficient or infeasible. Engineering applications arise  in power systems \cite{Dorfler2017},  communication networks \cite{Boyd2004}, machine learning \cite{Boyd2010} and  robotics \cite{Rabbat}, just to name a few. 
A prominent role is played by asynchronous algorithms, where communication and updates of the local processors are not  coordinated. The asynchronous approach is advantageous in several ways: it eliminates the need for synchronization, which is costly in large networks; it reduces the idle time, when distinct processors have different computational capabilities; it enhances robustness with respect to unreliable, lossy and delayed, communication; and it alleviates transmission and memory-access congestion. 
On this account, in this paper we  study a completely asynchronous implementation of the distributed dual ascent, a fundamental algorithm for constrained optimization.

\emph{Literature review:} The dual ascent consists of  solving the dual problem via the gradient method.  Its major advantage with respect to augmented Lagrangian methods (e.g., method of multipliers and \gls{ADMM}) is decomposability: for separable problems, the update breaks down into decentralized  subproblems, allowing for distributed and parallel (as opposed to sequential) implementation. Although this ``dual decomposition'' is an old idea \cite{Everett1963}, distributed algorithms based on the dual ascent are still very actively researched \cite{Necoara2015}, \cite{Beck2014}, \cite{NedicOzdaglar2009}.

Block-coordinate versions of the dual ascent, where only  part of the  variables is updated at each iteration, are  also explored in the literature \cite{Ananduta2020}. \blue{More generally, a variety of distributed algorithms has been proposed to solve constraint-coupled optimization problems, possibly with block-updates and time-varying communication \cite{CamisaCDC2019,FalsonePrandiniLCSS2021,Li_etal_TAC2021_Proxconstrainedopt}.  Nonetheless, in all the cited works, a common clock is employed to synchronize the communication and update frequencies. }

On the contrary, a global clock is superfluous for asynchronous methods. Since the seminal work  \cite{Bertsekas1989parallel}, distributed asynchronous optimization algorithms have received increasing attention, especially with regards to primal schemes \cite{SrivastavaNedic2011}, \cite{TianSunScutari2020}.
Most of the available asynchronous {dual} approaches leverage randomized activation \cite{Notarnicola2018}, \cite{Bastianello2020}, \cite{Lin_etal_2019}, \cite{Farina2019},
 where the update of each agent is triggered by a local timer or by signals received from the neighbors. However, no delay is tolerated, i.e., the agents perform their computations using the most recent information. This requires some coordination, since each agent must wait for its neighbors to complete their tasks before starting a new update. 
 To deal with delays caused by imperfect communication or non-negligible computational time, primal-dual \cite{Wu2018} and \gls{ADMM}-type \cite{Chang2016}, \cite{ARock} algorithms have been analyzed. Yet, the method in \cite{Chang2016} requires the presence of a master node; meanwhile, in  \cite{ARock}, \cite{Wu2018}, the delays are assumed to be independent of the activation sequences, which is  not realistic \cite{Wu2018}, \cite{TianSunScutari2020}.  
 
 \emph{Contributions:} In this paper, we propose an asynchronous implementation of the distributed dual ascent, according to the celebrated \emph{partially asynchronous} model, devised by Bertsekas and Tsitsiklis   \cite[§7]{Bertsekas1989parallel}, where: (a) agents perform updates and send data at any time, without any need for coordination signals; (b) the agents use outdated information from their neighbors to perform their updates.
In particular,
 we allow some agents to transmit  more frequently,  process faster and execute more iterations than others, based exclusively on their local clocks.
%
%
The model also encompasses the presence of heterogeneous delays or dropouts in the communication. 
Differently from  \cite{Chang2016}, only peer-to-peer communication is required.  Moreover, we do not postulate a stochastic characterization of  delays or activation sequences.
Instead, we merely assume bounds on the communication and update frequencies. 
In this partially asynchronous scenario, Low and Lapsley \cite{Low1999} studied a dual ascent algorithm, 
for a \gls{NUM} problem, with affine inequality constraints modeling link capacity limits. 

Here we consider a more general setting, and address the presence of convex  inequality and equality constraints.
Our main contribution is to  prove the convergence of the sequence generated by the  asynchronous distributed dual ascent to an optimal primal solution, under  assumptions that are standard for its synchronous counterpart and
provided that small-enough uncoordinated step sizes are employed. In fact, we drop some of the technical conditions in \cite{Low1999} (see §\ref{sec:problemstatement}), and relax the need for a global step. Our strategy is to relate the algorithm to a perturbed projected  scaled gradient scheme. 
\blue{With respect to  asynchronous primal methods \cite[§7.5]{Bertsekas1989parallel} and to the general fixed-point algorithm in \cite{Hannah2018}, the main  technical challenge is that
 the agents are not able to compute the partial gradients of the dual function locally; as a consequence, we have to consider
 two layers of delays}. To validate our results, we provide a numerical simulation on the  \gls{OPF} problem.

\emph{Notation}: 
$\mathbb{N}$ is the set of natural numbers, including $0$. 
$\overline{\R} := \R \cup \{\infty\}$ is the set of extended real numbers. 
$\0_n$ $\in \R^n$  is the vector with all elements equal to $0$; $I_n\in\R^{n\times n}$ is an identity matrix; we may omit the subscripts if there is no ambiguity. 
For an extended value function $f:\R^n \rightarrow \overline{\R} $, $\dom(f):=\{x\in\R^n \mid f(x)<\infty \}$; $f$ is  $\rho$-strongly convex if $x\mapsto f(x)-\frac{\rho}{2}\|x\|^2$ is convex, coercive if $\lim_{\|x\|\rightarrow \infty} f(x)=\infty$.  
Given a positive definite matrix $P\succ 0$, 
$\|{x}\|_{P}=:\sqrt{\langle x, P x \rangle}$ is the $P$-weighted norm;
we omit the subscript  if $P=I$. 
$\operatorname{int}(S)$ is the  interior of a set $S$.
\section{Problem setup and mathematical  background}\label{sec:problemstatement}
Let $\mc{I}=\{1,2,\dots,N\}$ be a set of agents. The agents can communicate over an undirected network $\mc{G}$($\mc{I},\mc{E})$; the  pair $(i,j)$ belongs to the set of edges $\mc{E}\subseteq \mc{I}\times \mc{I}$ if and only if agents $i$ and $j$ can \emph{occasionally} exchange information, with the convention $(i,i)\in\mc{E}$ for all $\i$. We denote by $\neigi=\{j \mid (i,j)\in\mc{E}\}$ the neighbors set  of agent $i$. The agents' common goal is to solve the following convex monotropic  optimization problem, where the decision vector  $x_i\in\R^{n_i}$ of agent $i$ is coupled to the decisions of the neighbors $\mc{N}_i$ via convex shared constraints:
\begin{subequations}\label{eq:problem}
	\begin{align}
		\label{eq:problem:minimization}
		\min_{ x_i\in\R^{n_i},{\i}  }  ~~& \sum_{\i} f_i(x_i) \\  	\textstyle 
		\text{s.t.}  ~~~&    \sum_{j\in\neigi} c_{i,j}(x_j) \leq \0_{p_i},~~~\forall \i \\ 	 
		 & \sum_{j\in\neigi} a_{i,j}(x_j)=\0_{r_i},~~~\forall \i, 
		\label{eq:problem:inequality} 
	\end{align}
\end{subequations}
Here, the cost $f_i:\R^{n_i}\rightarrow \overline{\R}$, the functions $\{ c_{j,i}:\R^{n_i}\rightarrow \R^{p_j}, j\in\mc{N}_i \}$, and the affine functions $\{ a_{j,i}:\R^{n_i}\rightarrow \R^{r_j}:x_i\mapsto A_{j,i}x_i-b_{j,i} , j\in\mc{N}_i \}$ are local data kept by agent $i$.
\begin{remark} \label{rem:applications}
	Problems in the form \eqref{eq:problem} arise naturally in resource allocation  \cite{Boyd2004}  and network flow problems \cite{Rockafellar1998}, e.g., \gls{NUM} for communication networks \cite{Beck2014} or \gls{OPF} in energy systems \cite{Dorfler2017}.  More generally, a well-known approach to solve a  (cost-coupled) distributed optimization problem is to recast it as \eqref{eq:problem}, by introducing slack variables to decouple the costs and additional consensus  constraints for consistency. For instance, the problem
	$\{ \min_{z\in \R} \sum_{\i} f_i(z) \}$ is equivalent to 
	\begin{equation}\label{eq:reformulation}
	\min_{x_i\in \R, \i} \ \textstyle \sum_{\i} f_i(x_i)  \text{ s.t. } L_{-1}(\col ((x_i)_{\i}))=\0_{N-1}  ,	\end{equation} 
	where $L_{-1}$ is the  (full row rank, see \cref{asm:fullrowrank} below) matrix obtained by removing the first row from the Laplacian of a connected graph; indeed, \eqref{eq:reformulation} is an instance of \eqref{eq:problem}. \hfill $\square$
\end{remark}

In the following, we use the compact notation $x:=\col((x_i)_{\i})\in\R^n$ and $f(x):=\textstyle \sum_\i f_i(x_i)$, where $n:=\sum_{\i} n_i$. Let us also define $c_{i,j}(x_j) := \0_{p_i}$, $A_{i,j}:=\0_{r_i \times n_j}$ and $a_{i,j}(x_j):=\0_{r_i}$  for all $\i$, $j\notin \mc{N}_i$; $m_i:=p_i+r_i$, $g_{i,j}(x_j):=\col(c_{i,j}(x_j), a_{i,j}(x_j)) $, $g_i(x_i):=\col((g_{j,i}(x_i))_{j\in\mc{N}})$, and $\Omega_i:=\R_{\geq 0}^{p_i} \times \R^{r_i}$, for all $i,j \in\mc{I}$; $g(x):=\sum_{\i} g_i(x_i)$, $m:=\sum_{\i} m_i$ and $\Omega:=\prod_{\i} \Omega_i$. 
We assume the following conditions throughout the paper. 

\begin{assumption}[Regularity and Convexity]\label{asm:blanket}~
	\begin{asmlist}
		\item \label{asm:convexity} For all $i\in\mc{I}$, $f_i$ is proper, closed,  and $\rho_i$-strongly convex, for some $\rho_i>0$.  
		 	\item \label{asm:smoothness} For all $i,j\in \mc{I}$,  $g_{i,j}$ is componentwise convex and $\theta_{i,j}$-Lipschitz continuous on $\dom(f_j)$, for some $\theta_{i,j}>0$.
		\item \label{asm:Slater}
			There exists  $\bar{x}\in\operatorname{int}(\dom(f))$ such that, for all $\i$,   $\sum_{j\in\neigi} c_{i,j}(\bar{x}_j)<  \0_{p_i}$ and $\sum_{j\in\neigi} a_{i,j}(\bar{x}_j)=\0_{r_i}$.
			\item \label{asm:fullrowrank}
			The matrix $A:=[A_{i,j}]_{i,j\in\I}$ has full row rank. \hfill $\square$
		\end{asmlist}
\end{assumption}
Under \cref{asm:Slater}, problem \eqref{eq:problem} is feasible.  In addition, the strong convexity in  \cref{asm:convexity} ensures that  there exists a unique solution $x^\star$, with finite optimal value $f^\star:=f(x^\star)\in\R $; this condition is standard for dual gradient methods \cite[Asm.~2.1]{Necoara2015}, \cite[Asm.~1]{Ananduta2020}, and commonly found in the problems mentioned in \cref{rem:applications}. Differently from \cite[Asm.~C1]{Low1999}, we do not assume  differentiability of $f$, nor that the functions $f_i$'s are increasing (e.g., quadratic cost functions are allowed here). We note that 
 local constraints can be  enforced in \eqref{eq:problem} by opportunely choosing the domains  of the $f_i$'s (which need not  be bounded, cf. \cite{Low1999}). We also remark that \cref{asm:smoothness} is  automatically satisfied in the most common case of affine inequality constraints \cite{Beck2014}, \cite{Necoara2015}. 

Given the information available to each agent, the  natural way of distributedly solving  \eqref{eq:problem} is resorting to dual  methods. 
By \cref{asm:Slater}, strong duality holds \cite[Th.~28.2]{Rockafellar1970ConvexAnalysis}, i.e.,
\begin{align}\label{eq:strongduality}
f^\star=q^\star:=\underset{y\in \Omega}{\max}  \ q(y), 
\end{align}
where $q: \Omega \rightarrow \R $ is the concave dual function, 
\begin{equation}\label{eq:dualfunction}
q(y):=\underset{x\in\R^n}{\min} \ f(x)+\langle g(x), y \rangle,\end{equation}
with dual variable $y$, and the maximum $q^\star$ is attained in \eqref{eq:strongduality}; we denote by $\mc{Y}^\star$ the convex nonempty set of dual solutions (namely, solutions of \eqref{eq:strongduality}). Moreover,
\cref{asm:fullrowrank} rules out the case of redundant equality constraints; together with \cref{asm:Slater}, it guarantees that the (convex) function $-q$ is coercive on $\Omega$, and hence that $\mc{Y^\star}$ is  bounded (for similar arguments, see \cite[§VII, Th.~2.3.2]{Hiriart-Urruty2012}, \cite[Lem.~1]{NedicOzdaglar2009}).
 For this reason, \cref{asm:Slater,asm:fullrowrank} have been exploited for particular instances  of problem~\eqref{eq:problem}  \cite[Asm.~2.1]{Necoara2015}, \cite[Asm.1]{NedicOzdaglar2009}.
%
\subsection{Synchronous distributed dual ascent algorithm}
 Under \cref{asm:convexity}-(ii), the concave dual function $q$ in \eqref{eq:dualfunction} is differentiable with Lipschitz  gradient {\cite[Lem. II.1]{Beck2014}}. Thus, for a  $\gamma>0$ small enough, the dual ascent iteration
\begin{align}\label{eq:synchronousda}
y(k+1)=\proj_{\Omega} \left ( y(k)+\gamma \nabla q (y(k))\right )
\end{align}
converges to a dual solution. By the envelope theorem, 
$\nabla q(y)=g(x^\star(y))$, where $x^\star(y):=\argmin_{x\in\R^n} \  f(x)+\langle g(x), y \rangle$; therefore, by assigning to agent $i$ the Lagrangian multipliers $y_i\in\R^{m_i}$, with $y=\col((y_i)_{\i})$, the update in \eqref{eq:synchronousda} can be written, for the single agent $i$, for all $\i$, as
\begin{subequations}\label{eq:synchronousdai}
	\begin{align}
	\label{eq:synchronousdai1}
x_i(k+1) &=   \underset{ x_i\in \R^{n_i}}{\argmin} \, \Big(
f_i(x_i)+  \underset{j\in \neigi }{\textstyle \sum} \langle g_{j,i}(x_i), y_j(k) \rangle  \Bigr) \hspace{-1em}
\\[-0.5em]
\label{eq:synchronousdai2}
y_i(k+1)&=\proj_{\Omega_i} \,\Bigl(  y_i(k)+ \gamma   \underset{j\in \neigi }{\textstyle \sum}  g_{i,j} (x_j(k+1))   \Bigr),
\end{align}
where the $\argmin$ is single valued by \cref{asm:convexity}, and the sequence $(x(k))_\k$ converges to the primal solution $x^\star$. We emphasize that computing the update in \eqref{eq:synchronousdai} requires each agent to receive information \emph{from all of its neighbors  twice} per iteration: the first time  because computing $x_i(k+1)$ requires the knowledge of  $\{y_j(k),{j\in\mc{N}_i}\}$; and the second time to collect the   quantities $\{g_{i,j}(x_j(k+1)),{j\in\mc{N}_i}\}$.
\end{subequations}

\section{Asynchronous distributed   dual ascent}
\floatstyle{spaceruled}
\restylefloat{algorithm}
\begin{algorithm*}[h]
	\caption{Asynchronous Distributed Dual Ascent} \label{algo:2}
	\vspace{0.5em}
	\textbf{Initialization:} $\forall \i$, $y_i(0)=\0_{m_i}$, $x_i(0)=\underset{x_i\in \R^{n_i}}{\argmin} \ f_i(x_i)$. \\
	
	\textbf{Local variables update:} For all $\k$, each agent $\i$ does:
\begin{subequations}
	\label{eq:algo:2_TI}
\begin{empheq}[left={\text{if } k\in T_i: \hspace{5em}\empheqlbrace \  }]{align}
		\label{algo:xupdate}
x_i(k+1) &=\underset{x_i\in \R^{n_i}}{\argmin} \ \Biggl(
f_i(x_i)+  \sum_{j\in \neigi } \langle g_{j,i} (x_i), y_j(\tau_{i,j}^k)  \rangle \Biggr)
\\
\label{algo:yupdate}
y_i(k+1) &={\proj}_{\Omega_i} \Biggl( \ y_i(k)+ \gamma_i   \sum_{j\in \neigi } g_{i,j} (x_j(\tau_{i,j}^k) )  \Biggr) \hspace{10em}
	\end{empheq}
\end{subequations}
\vspace{-0.5em}
\begin{subequations}
	\label{eq:algo:2_notTI}
	\begin{empheq}[left={\text{otherwise}:  \hspace{5em}\empheqlbrace \  }]{align}
		x_i(k+1)& =x_i(k) \label{algo:xupdate2}\\
y_i(k+1) &=y_i(k). \hspace{24.5em} \label{algo:yupdate2}
	\end{empheq}
\end{subequations}
\vspace{-0.5em}
\end{algorithm*}

Let us now introduce an  asynchronous version of the distributed dual ascent method \eqref{eq:synchronousdai}, according to the 
\emph{partially asynchronous model} \cite[§7.1]{Bertsekas1989parallel}. 
 The iteration is shown in \cref{algo:2}, and it is  determined by: 
\begin{itemize}[leftmargin=*]
	\item a nonempty  sequence $T_i \subseteq \N$, for each $\i$. Agent $i$ performs an update only for $k\in T_i$.
	\item an integer variable $\tau_{i,j}^k $, with  $0 \leq \tau_{i,j}^k \leq k$, for each $\i$, $j\in\mc{N}_i$, $k\in \N$, which represents the number of steps by which the information used in the updates of $(x_i, y_i)$ at step $k$ is outdated. For example, $\tau_{i,j}^k=k-\delta$ means that the variable $y_j(\tau_{i,j}^k)$ that agent $i$ uses to compute $x_i({k+1})$ is outdated by $\delta$ steps. 
%
\end{itemize}
In particular, the following variables are available to agent $i$ when performing the update at $k\in T_i$: 
\begin{align}
	y_i(k), \  \{ y_j(\tau_{i,j}^k) \mid  j\in\neigi \}, \ 
	\{ g_{i,j} (x_j(\tau_{i,j}^k)) \mid  j\in \neigi \}.
\end{align}
%
%
%
%
For mathematical convenience, $\tau_{i,j}^k$ is defined also for $k\notin T_i$, even if these variables are  immaterial to \cref{algo:2}. 
Note that $\tau_{i,j}^k\geq 0$ implies that at  $k=0$ all the agents have  updated information on the neighbors' variables; this  is not restrictive and it is only meant to ease the notation (the same holds for the initialization in \cref{algo:2}; see also \cite[§7.1]{Bertsekas1989parallel}).
\begin{remark}\label{rem:alwaysdelay}
	The update in   \eqref{eq:algo:2_TI} differs from \eqref{eq:synchronousdai}  even if $\tau_{i,j}^k=k$ for all $\i$, $j\in \mc{N}_i$, $\k$, as in  \eqref{eq:synchronousdai2} the agents are exploiting the most recent information $x_j(k+1)$.
	\hfill $\square$ 
\end{remark}

\begin{figure}
\centering
\vspace{0.5em}
\includegraphics[width=1\columnwidth]{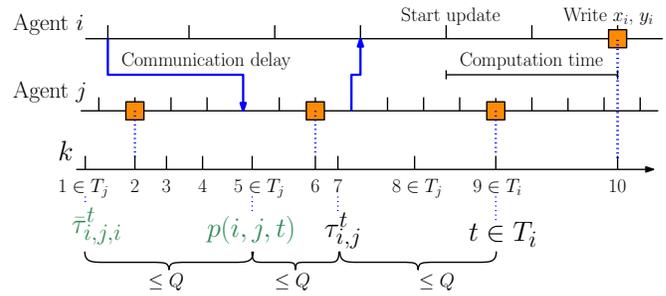}
\caption{An example of asynchronous updates and communication. We consider a simple scenario, where each agent periodically  broadcasts its variables and computes an update using the latest data received from its neighbors, according to its \emph{local} clock. The communication is subject to bounded delays.  The figure only considers two neighboring agents $i$ and $j$.   The  squares indicates the instants when the local variables change value; the arrows indicates  data transmission. The global counter $k$ increases every time an agent in the network completes its update. Here,  agent $i$  completes an update using the outdated information $x_j(7), y_j(7)$ and  $k$ is increased to $10$:  hence  agent $i$ is viewed by an external observer as performing the update at $k=9$, and $\tau_{i,j}^9=7$. \blue{The quantities in green are defined in §\ref{sec:convergence}.}
}
\label{fig:asyn_up}
\end{figure}

In \cref{algo:2}, $k$ should be regarded as an event counter, an iteration index as seen by an external observer, which is increased every time one or more agents complete an update. It is introduced to give a global  description of the algorithm, but it is \emph{not known} by the agents, who can compute and communicate \emph{without any form of coordination}. We present a simple example in \cref{fig:asyn_up}.  Indeed, the  partially asynchronous model  captures a plethora of asynchronous  protocols (for different choices of the parameters $T_i$'s, $\tau_{i,j}^k$'s), encompassing several scenarios where: (a)~the communication links of the network $\mc{G}$ are lossy or active intermittently; (b)~there are heterogeneous delays on the transmission; (c)~the local computation time cannot be neglected; (d)~the agents update their local variables at different frequencies. We refer to \cite{Bertsekas1989parallel} for an exhaustive discussion.
\begin{assumption}[Partial asynchronism, {\cite[§7.1, Asm. 1.1]{Bertsekas1989parallel}}] \label{asm:partial asynchronicity}
	There exists a positive integer $Q$ such that:
	\begin{itemize}[topsep=0pt]
		\item[(i)] \emph{Bounded inter-update intervals}: for all $\k$, for all $\i$, it holds that $\{k,k+1,\dots,k+Q-1\} \cap T_i \neq \varnothing$;
		\item[(ii)] \emph{Bounded  delays}: for all $\k$, for all $i \in\I$ \blue{and for all $j\in\neigi$}, it holds that  $k-Q \leq \tau_{i,j}^k \leq k$. \hfill $\square$
	\end{itemize}
\end{assumption}
\begin{remark}
	\cref{asm:partial asynchronicity} is mild and easily satisfied in distributes computation; for instance, it holds for the  scenario in \cref{fig:asyn_up}, see   \cite[§7.1, Ex~1.1]{Bertsekas1989parallel},
	\cite[§III.A]{TianSunScutari2020}. \hfill $\square$
\end{remark}

We are ready to enunciate our main result. 
\begin{theorem}\label{th:main}
	Let Assumptions \ref{asm:blanket}-\ref{asm:partial asynchronicity} hold. For all $\i$, let 
	\begin{align}
	\label{eq:const_thetai}
	\theta_i &:=\sqrt{\textstyle \sum_{j \in \mc{N}_i } \theta_{j,i }^2}\\
	\label{eq:const_phii}
	\phi_i&:=\textstyle\sum_{j\in\mc{N}_i} \frac{\theta_j ^2}{\rho_{\blue{i}}}\\
	\ell_i&:=\textstyle\sum_{j\in\neigi} \theta_{i,j} \frac{\theta_j}{\rho_j} \label{eq:const_elli}\\
	\label{eq:const_xii}
	\xi_i&:=\blue{\textstyle\sum_{j\in \neigi } \sum_{l\in\mc{N}_j} \theta_{l,j}{\textstyle\frac{\theta_j}{\rho_j}}}
	\end{align}
	with $\theta_{j,i}$ and $\rho_j$ as in \cref{asm:blanket}, for all $i,j\in\I$. Assume that, for all $\i$, the step size $\gamma_i>0$ is chosen such that
	\begin{align}\label{eq:gamma_bound}
	\textstyle\gamma_i^{-1} > \frac{1}{2}\phi_i+\blue{\frac{3}{2}Q (\ell_i+\xi_i)},
	\end{align}
	with $Q$  as in \cref{asm:partial asynchronicity}. Then, the sequence $(x(k))_{\k}$ generated by \cref{algo:2} converges to the solution $x^
	\star$ of the optimization problem in  \eqref{eq:problem}. \hfill $\square$
\end{theorem}

\begin{remark}
	Each agent $i$ can locally compute the parameters $\phi_i$, $\ell_i$, $\xi_i$ in \cref{th:main}, only based on some information  from its direct neighbors. Thus, the choice of the step sizes $\gamma_i$'s is decentralized, provided that the agents have access to (an upper bound for) the asynchronism bound $Q$ (otherwise, vanishing step sizes can be considered). \color{black} \hfill $\square$
\end{remark}
\begin{remark}
	 \blue{
	As usual for partially asynchronous optimization algorithms \cite{Bertsekas1989parallel}, \cite{TianSunScutari2020}, the upper bound  in \eqref{eq:gamma_bound} decreases if $Q$ grows. This is a structural issue: we can construct a problem satisfying \cref{asm:blanket} (similarly to \cite[§7.1, Ex.1.3]{Bertsekas1989parallel}) such that, for any fixed  positive $\gamma_i$'s, there is a large enough $Q$ and some sequences $\tau_{i,j}^k$'s , $T_i$'s satisfying \cref{asm:partial asynchronicity}, for which \cref{algo:2} diverges.  
  \hfill $\square$}
\end{remark}

\section{Convergence analysis}\label{sec:convergence}
In this section we prove \cref{th:main}. Our idea is to relate \cref{algo:2} to a perturbed scaled block coordinate version of \eqref{eq:synchronousda}, and to show that, for small-enough steps sizes, the error caused by the outdated information is also small and does not compromise convergence.
\blue{With respect to the asynchronous gradient method in \cite{Bertsekas1989parallel},
the main technical complication is that, for each update in \eqref{eq:synchronousda}, the agents   communicate twice;  in turn, the update of agent $i$ depends  on the variables of its neighbors' neighbors (or \emph{second order neighbors}). In the following, we denote by  $\tsup[1]{y}_i:=\col((y_j)_{j \in \mc{N}_i})$  and $\tsup{y}_i:=\col((\tsup[1]{y}_j)_{j\in \mc{N}_i})$ the dual variables of the neighbors and second order neighbors of agent $i$, respectively. }

\blue{To compare \cref{algo:2} and \eqref{eq:synchronousda}, the first step is to get rid of the primal variables $x_i$'s in \cref{algo:2}.
Importantly, we have to take into  account that the $x_j$'s in \eqref{algo:yupdate} are not only outdated, but also computed using outdated information. 
In fact,} for any $i\in \mc{I}$, $k\in T_i$,
%
%
%
%
%
%
%
%
%
%
%
%
the variable $x_j(\tau_{i,j}^k)$  is computed by agent $j \in \mc N_i$, according to \eqref{algo:xupdate}, as
\begin{align}
\label{eq:xitau}
x_j(\tau_{i,j}^k)& =\hphantom{:}\underset{x_j\in \R^{n_j}}{\argmin} \ 
f_j(x_j)+  \sum_{l\in \mc{N}_j } \langle g_{l,j}( x_j), y_l(\bar{\tau}^k_{i,j,l})  \rangle  \\
&=: \h_j(\tsup[1]{y}_j (\bar{\tau}^k_{i,j})), \label{eq:hj}
\end{align}
where, for all $l\in\mc{N}_j$,   \blue{$\bar{\tau}^k_{i,j,l}:={\tau}^{p}_{j,l}$, and $p$ is the last time agent $j$ performed an update prior to  $\tau_{i,j}^k$ (see \cref{fig:asyn_up}),} i.e., 
\begin{align}\label{eq:bartau}
 \hspace{-1em} p= p(i,j,k):=\max ( \{0\} \cup  \{t \in T_j \mid t \leq  \tau_{i,j}^k-1\} ),
\end{align}
 (and  \eqref{eq:xitau}  also holds if  $p(i,j,k)=0$ because of  the initial conditions in \cref{algo:2}). For brevity of notation, 
 in  \eqref{eq:hj}, we define 
 $\bar{\tau}^k_{i,j}:=\col((\bar{\tau}^k_{i,j,l})_{l\in\mc{N}_j})$ and $\tsup[1]{y}_j(\bar{\tau}^k_{i,j}):= \col( (y_l(\bar{\tau}^k_{i,j,l})_{l\in\mc{N}_j})$.
 \blue{By replacing  \eqref{eq:hj} in \eqref{algo:xupdate},} we obtain that,   for all $\i$,  $k\in T_i$
\begin{align}
\nonumber
y_i(k+1)&= \hphantom{:} {\proj}_{\Omega_i} \Biggl( \ y_i(k)+ \gamma_i  \sum_{j\in \neigi } g_{i,j}(\h_j(\tsup[1]{y}_j (\bar{\tau}^k_{i,j}))) \Biggr)\\
&=:{\proj}_{\Omega_i} \left( \ y_i(k) + \gamma_i 
 \bs{F}_i(\tsup{y}_i(\bar{\tau}^k_{i}) ) \right), \label{eq:compactalgo1}
\end{align}
where $\bar{\tau}_i^k:=\col((\bar{\tau}^k_{i,j})_{j\in\mc{N}_i})$,  
$\tsup{y}_i(\bar{\tau}_i^k):=\col ((\tsup[1]{y}_j(\bar{\tau}^k_{i,j}))_{j\in\neigi})$. \blue{Hence, \eqref{eq:compactalgo1} expresses the update in \eqref{algo:yupdate} as a function of the outdated dual variables  of the second order neighbors of agent $i$.}
Since \eqref{eq:compactalgo1} makes use of second order information,  the maximum delay (in this notation) is not bounded  by $Q$ anymore: instead, by  \eqref{eq:bartau} and \cref{asm:partial asynchronicity}, \blue{it holds that $p(i,j,k) -Q \leq \tau_{j,l}^{p(i,j,k)}= \bar{\tau}_{i, j, l}^{k} \leq p(i,j,k) $, $\tau^k_{i,j}-Q \leq  p(i,j,k) \leq \max (0, \tau^k_{i,j}-1) $, $k-Q \leq  \tau^k_{i,j} \leq k$;} thus,
 \begin{equation}\label{eq:alwaysdelay}
 k-3Q\leq \bar{\tau}^k_{i,j,l} \leq k-1,
 \end{equation}
 for all $k\geq 1$, and $\bar{\tau}^0_{i,j,l}=0$.
 Figure \ref{fig:asyn_up} illustrates how the lower bound is obtained, with $l=i$.
 
 	\blue{We emphasize that the mapping $\bs{F}_i$ is not a partial gradient of the dual function (differently from the asynchronous gradient method in  \cite{Bertsekas1989parallel}) -- e.g., its argument lies in a different space.
 In particular, we note that   $\tsup{y}_i(\bar{\tau}^k_{i})$ can contain multiple instances of $y_l$, for some $l\in\mc{I}$ (including $y_i$), with different delays. Nonetheless,}
by the definition in \eqref{eq:compactalgo1},  for any $y(k)\in\R^m$, we have
\begin{equation}\label{eq:onconsensus}
\bs{F}_i(\tsup{y}_i(k))=\nabla_{y_i} q(y(k)).
\end{equation}
In one word, if there is no delay, \eqref{eq:compactalgo1} corresponds to the $y_i$-update of the synchronous dual ascent \eqref{eq:synchronousda} (however, there is always delay in \cref{algo:2}, see \eqref{eq:alwaysdelay} and \cref{rem:alwaysdelay}). 

Finally, we can also rewrite \cref{algo:2} as 
\begin{equation}\label{eq:compactform}
	(\forall \k)(\forall \i) \quad y_i(k+1)=y_i(k)+\gamma_i s_i(k),
\end{equation}
where 
\begin{equation}\label{eq:si}
s_i(k):= \frac{1}{\gamma_i} ( {\proj}_{\Omega_i} ( \ y_i(k) + \gamma_i
\bs{F}_i(\tsup{y}_i(\bar{\tau}^k_{i}) )  ) -y_i(k)),
\end{equation}
if $k\in T_i$,  $s_i(k)=\0_{m_i}$ otherwise. 

Before proceeding with the analysis of \cref{th:main}, we recall the following results. The proof of \cref{lem:hlip,} is standard and  omitted here (see, e.g., \cite[Lem.~1]{Ananduta2020}).
\begin{lemma}\label{lem:hlip}
	For all $j \in \mc{I}$, the mapping $\h_j$  in \eqref{eq:hj} is $\frac{\theta_j}{\rho_j}$-Lipschitz continuous, with $\theta_j$ as in \eqref{eq:const_thetai}.   \hfill $\square$
\end{lemma}
\begin{lemma}[Weighted descent lemma, {\cite[Lem.~2.2]{Necoara2015}}]
\label{lem:descentlemmaPhi}
	Let  $\Phi=\diag((\phi_i \otimes I_{m_i})_\i)$, 
	 $\phi_i$ as in \eqref{eq:const_phii},
	for all $\i$.  Let $q$ be the dual function in  \eqref{eq:dualfunction}. For any $y,z\in\R^m$, it holds that
 \[  q(y) \geq  q(z)+\langle  y-z, \nabla q (z ) \rangle -
  \textstyle	\frac{1}{2} \|y-z\|^2_\Phi.    \QEDopenhereeqn
 \]
\end{lemma}
\smallskip
\begin{lemma}[{\cite[Lemma 5.1]{Bertsekas1989parallel}}]\label{lemma:siconsistency}
	Let  $s_i(k)$ be as in \eqref{eq:si}, for all $\i$. Then, for all $\i$, $\k$, it holds that 
	\[
	\langle s_i(k), \bs{F}_i(\tsup{y}_i(\bar{\tau}^k_{i}) )\rangle  \geq  \|s_i(k)\|^2.
	\QEDopenhereeqn
	\]
\end{lemma}

\subsection{Proof of \cref{th:main}}
	\noindent Let $\Gamma:= \diag( (\gamma_i \otimes I_{m_i})_{\i})$, $s(k)=\col((s_i(k))_\k)$, $\bs{F}(\tsup{y}(\bar{\tau}^k)):=\col(    (\bs{F}_i(\tsup{y}_i(\bar{\tau}^k_{i}) ))_{\i})$.
	By \cref{lem:descentlemmaPhi}, we have 
	\begin{align*}
	 q(y(k+1)) & =q(y(k)+\Gamma s(k)) \\
	& \geq   q(y(k)) -	\textstyle \frac{1}{2} \|\Gamma s(k)\|^2_\Phi
	+\langle  \Gamma s(k), \nabla q (y(k)) \rangle \\ 
	& = q(y(k))-	\textstyle\frac{1}{2} \|\Gamma s(k)\|^2_\Phi
	+\langle  \Gamma s(k),   \bs{F}(\tsup{y}(\bar{\tau}^k) )\rangle \\
	& \qquad +\langle  \Gamma s(k), \nabla q (y(k)) - \bs{F}(\tsup{y}(\bar{\tau}^k) )\rangle  \\
	& = q(y(k))- \textstyle\frac{1}{2} \|\Gamma  s(k)\|^2_\Phi
	+\langle  \Gamma s(k),  \bs{F}(\tsup{y}(\bar{\tau}^k) )\rangle \\
	& \qquad +\langle  \Gamma s(k), \bs{F}(\tsup{y}(k)) -\bs{F}(\tsup{y}(\bar{\tau}^k) )\rangle  \\
	& \geq   q(y(k))- \textstyle\frac{1}{2} \| s(k)\|^2_{\Phi\Gamma^2}+
	\| s(k)\|_{ \Gamma}^2 \\
	\numberthis 
	& \qquad +\langle  \Gamma s(k), \bs{F}(\tsup{y}(k)) -\bs{F}(\tsup{y}(\bar{\tau}^k) )\rangle,
	\hspace{-1.5em} \label{eq:usefulstep}
 	\end{align*} \color{black}
where in the last equality we used \eqref{eq:onconsensus} and the last inequality follows by \cref{lemma:siconsistency} (we  recall that $\Phi, \Gamma$ are diagonal matrices).
 We next bound the last addend in \eqref{eq:usefulstep}.  
\blue{ By definition of $ \bs{F}_i$ in \eqref{eq:compactalgo1} and the Cauchy--Schwartz inequality, we have
 \begin{align*}
&\langle  \Gamma s(k), \bs{F}(\tsup{y}(k)) -\bs{F}(\tsup{y}(\bar{\tau}^k) )\rangle 
\\ &
 \geq - \sum_{\i}  \| \gamma_i s_i(k) \| \| \bs{F}_i(\tsup{y}_i(k))-\bs{F}_i(\tsup{y}_i(\bar{\tau}^k_{i})) \|  
 \\ &
\geq -\sum_{\i}  \| \gamma_i s_i(k) \| \! \sum_{j\in \neigi } \! \! \|  g_{i,j}(\h_j(\tsup[1]{y}_j (k))) \! - \! g_{i,j}(\h_j(\tsup[1]{y}_j (\bar{\tau}^k_{i,j}))) \|
\\ &
 \geq -\sum_{\i}  \| \gamma_i s_i(k) \| \sum_{j\in \neigi } \theta_{i,j}{\textstyle\frac{\theta_j}{\rho_j}} \ \|\tsup[1]{y}_j (k)-\tsup[1]{y}_j (\bar{\tau}^k_{i,j}) \|
 \\&
 \geq -\sum_{\i}  \| \gamma_i s_i(k) \| \sum_{j\in \neigi } \theta_{i,j}{\textstyle\frac{\theta_j}{\rho_j}} \sum_{\tau=k-3Q}^{k-1} \!\! \| \col((\gamma_l s_l(\tau))_{l\in \mc{N}_j}) \|,
   \end{align*}
   where in the third inequality we used \cref{lem:hlip} and \cref{asm:smoothness}, and the last  follows by \eqref{eq:compactform} and \eqref{eq:alwaysdelay} (without loss of generality, we let $s(k):= \0_m$, for all $k<0$). Therefore,  by the elementary relation $2|a||b|\leq a^2+b^2$, we obtain
   \begin{align*}
 & \langle  \Gamma s(k), \bs{F}(\tsup{y}(k)) -\bs{F}(\tsup{y}(\bar{\tau}^k) )\rangle \\
& 
 \quad \geq  - \sum_{\tau=k-3Q}^{k-1}\sum_{\i} \sum_{j\in \neigi } \theta_{i,j}{\textstyle\frac{\theta_j}{\rho_j}} \   {\textstyle\frac{1}{2}} \| \gamma_i s_i(k) \|^2 \\
 & \numberthis \label{eq:twotermsstep}   \hphantom{ {}\leq{} } \quad  -\sum_{\tau=k-3Q}^{k-1}\sum_{\i} \sum_{j\in \neigi } \theta_{i,j}{\textstyle\frac{\theta_j}{\rho_j}}   \sum_{l\in\mc{N}_j} {\textstyle\frac{1}{2}}\| \gamma_l s_l(\tau) \|^2.
 \end{align*}
 The first term on the right-hand side of \eqref{eq:twotermsstep} equals $ -\|s(k) \|^2_{ \frac{3}{2}Q\Gamma^2 L} $, with $L := \diag((\ell_{i} \otimes I_{m_{i}})_{i \in \mathcal{I}})$ and $\ell_i$ as in \eqref{eq:const_elli}.
For the second term, since $\mc{G}$ is undirected, we  reorder the addends as
\begin{align*}
& 
 \sum_{\tau=k-3Q}^{k-1} \sum_{\i}  \sum_{j\in \neigi } \theta_{i,j}{\textstyle\frac{\theta_j}{\rho_j}}
\sum_{l\in\mc{N}_j} {\textstyle \frac{1}{2}}
\|\gamma_l s_l(\tau) \|^2
\\ & \quad = \sum_{\tau=k-3Q}^{k-1} \sum_{\i} {\textstyle \frac{1}{2}} \| \gamma_i s_i(\tau)\|^2 \left( \sum_{j\in \neigi } \sum_{l\in\mc{N}_j} \theta_{l,j}{\textstyle\frac{\theta_j}{\rho_j}} \right)
\\ & \quad = \sum_{\tau=k-3Q}^{k-1} \| s(\tau)\|^2_{{ \frac{1}{2}}\Gamma^2 \Xi},
\end{align*}
where $\Xi=:\diag((\xi_i \otimes I_{m_i})_\i)$ and $\xi_i$ as in \eqref{eq:const_xii}.}
Then, by substituting  in \eqref{eq:usefulstep}, we obtain 
\begin{align*}
q(y(k+1)) &\geq   q(y(k))+\| s(k)\|^2_{ \Gamma (I-\frac{1}{2}\Phi\Gamma-\frac{3}{2}Q L \Gamma)} \\ 
& \qquad - \sum_{\tau=k-3Q}^{k-1} \|s(\tau)\|^2_{{\textstyle \frac{1}{2}}\Xi \Gamma^2}  \numberthis \label{eq:usefulstep2},
\end{align*}
with  $I-\frac{1}{2}\Phi\Gamma-\frac{3}{2}Q L \Gamma \succ 0$ by the assumption on $\Gamma$ in \eqref{eq:gamma_bound}. Since \eqref{eq:usefulstep2} holds for any $k$, summing over $k$ finally yields
\begin{align*}
& q(y(k+1)) -q(y(0)) \\
& \qquad  \geq  \sum_{\tau=0}^{k} \| s(\tau)\|^2_{ \Gamma (I-\frac{1}{2}\Phi\Gamma-\frac{3}{2}Q L \Gamma)} - 3Q \sum_{\tau=0}^k  \|s(\tau)\|^2_{\frac{1}{2}\Xi \Gamma^2}
 \\ \numberthis
& \qquad = \sum_{\tau=0}^{k} \| s(\tau)\|^2_{ \Gamma (I-\frac{1}{2}\Phi\Gamma-\frac{3}{2}Q(\Xi+L)\Gamma)},  \label{eq:final inequality}
\end{align*}
and $ I-\frac{1}{2}\Phi\Gamma-\frac{3}{2}Q(\Xi+L)\Gamma\succ 0$ by  assumption  \eqref{eq:gamma_bound}. 

We know that $q$ is bounded above on $\Omega$ by \eqref{eq:strongduality}, and $y(k)\in\Omega\ $ for all $ \k$; then,  by taking the limit  in \eqref{eq:final inequality}, we have \begin{align}\label{eq:siconv}
\lim_{k\to\infty} s(k) = \0_m.
\end{align}
Hence, by the updates in \eqref{eq:compactform} and by \eqref{eq:siconv}, we also have 
\begin{equation}\label{eq:yresidualconvergence}
\lim_{k\to\infty} \| y(k+1)-y(k) \| =0.
\end{equation}
 Similarly, since, by \cref{asm:partial asynchronicity}(ii) and \eqref{eq:compactform}, $\| y_j(k)-y_j(\tau^k_{i,j}) \| \leq \sum_{\tau=k-Q}^{k-1} \gamma_j\|s_j(\tau)\|$, it also follows that 
\begin{equation}\label{eq:tauconv}
\lim_{k\to\infty} \| y_j(k)-y_j(\tau^k_{i,j}) \| =0, \ \forall \i, j\in \neigi.
\end{equation} 

For any $\i$, consider the \emph{subsequence} $(s_i(k))_{k\in T_i}$, which converges to $\0$ by \eqref{eq:siconv}. In view of \eqref{eq:si},  $( {\proj}_{\Omega_i} ( \ y_i(k)+ \gamma_i
\bs{F}_i(\tsup{y}_i(\bar{\tau}^k_{i}) )  ) -y_i(k))_{k\in T_i} \rightarrow \0$. Therefore, by leveraging the  continuity of $\bs{F}_i$ (which directly follows by the definition in \eqref{eq:compactalgo1} and \cref{lem:hlip}) and of the projection \cite[§3.3, Prop.~3.2]{Bertsekas1989parallel},  \eqref{eq:tauconv}  and  \eqref{eq:onconsensus} yield 
 $( {\proj}_{\Omega_i} ( \ y_i(k)+ \gamma_i
\nabla_{y_i}q(y(k)) -y_i(k)))_{k\in T_i }\rightarrow \0$. However, again by  continuity, \eqref{eq:yresidualconvergence} and \cref{asm:partial asynchronicity}(i), we can also infer convergence of the whole sequence,
$
\lim_{k\to\infty} \ {\proj}_{\Omega_i} ( \ y_i(k)+ \gamma_i
\nabla_{y_i}q(y(k)) -y_i(k)=\0_{m_i}
$, or 
\begin{equation}\label{eq:residualconv}
\lim_{k\to\infty} \ \left( {\proj}_{\Omega_i} ( \ y(k)+\Gamma
\nabla q(y(k)) -y(k) \right) =\0_{m}.
\end{equation}
We note that $q(y(k))\geq  q(y(0))$ for all $\k$, by \eqref{eq:final inequality}; moreover, $-q$ is coercive on $\Omega$ by \cref{asm:Slater,asm:fullrowrank}. We conclude that the sequence $(y(k))_\k$ is bounded; in turn, \eqref{eq:residualconv} implies that $(y(k))_\k$ converges to the set of dual solutions $\mc{Y}^\star$.


We can finally turn our attention to the primal problem \eqref{eq:problem}.
By \eqref{algo:xupdate}, for any $\i$, $k \in T_i$, we have $x_i(k+1)=\h_i(\tsup[1]{y}_i(\tau_i^k))$, where $\h_i$ is defined in \eqref{eq:hj} and $\tau_i^k:=\col((\tau_{i,j}^k)_{j\in\neigi})$. Moreover, by strong duality, for any $y^\star=\col((y_i^\star)_\i) \in \mc{Y^\star}$, it holds that  $\h_i(\tsup[1]{y}_i^\star) =x_i^\star$, with $\tsup[1]{y}_i^\star:= \col ((y_j^\star)_{j\in\neigi})$ and $\col((x_i^\star)_\i)=x^\star$. Therefore we can exploit  \eqref{eq:tauconv},
the fact that $(y(k))_k$ is converging to $\mc{Y}^\star$, and Lipschitz continuity of $\h_i$ in \cref{lem:hlip}, to conclude that $(x_i(k+1))_{k\in T_i} \rightarrow x_i^\star$. The conclusion follows because the convergence also holds for the whole sequence, i.e., 
$(x_i(k))_{\k} \rightarrow x_i^\star$, by \eqref{algo:xupdate2}.
\hfill$\blacksquare$ 

\section{Numerical simulation}
\begin{figure}
	\centering
	\vspace{1em}
	\includegraphics[width=1\columnwidth]{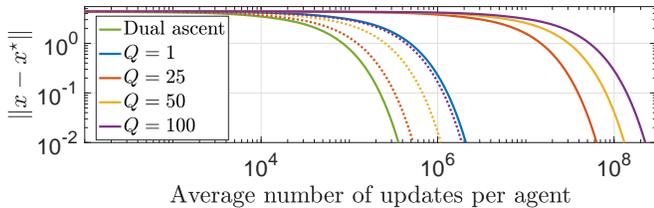}
	\caption{Distance from the optimum, \blue{with step sizes chosen  to satisfy their theoretical upper bounds (solid lines) and $100$ times larger (dotted lines).}
	}
	\label{fig:simulation}
\end{figure}
We consider  an \gls{OPF} problem on the IEEE 14-bus network \cite[Fig.~2]{Ananduta2020}. Each bus $i\in \mc{I}=\{1,\dots,14\}$ has a decision variable $x_i=\col(P_i,\psi_i)\in \R^2$, where $P_i\geq 0$ is the power generated, bounded by generation capacities, and $\psi_i$ is the  
voltage phase of bus $i$. The goal is to minimize the sum of strongly convex quadratic local costs, subject to the coupling flow constraints $\{P_i -P_i^\text{d} =\textstyle \sum_{j\in \mc{N}_i} B_{i,j} (\psi_i - \psi_j),  \forall \i \}$,
where $P_i^\text{d} \geq 0$ is the power demand at bus $i$ and $B_{i,j}$ is the susceptance of line $(i,j)$, which represent the  direct current approximation of the power flow equations. 
We simulate \cref{algo:2} for the setup in \cref{fig:asyn_up},
with randomly chosen delays and local clock frequencies. \blue{We compare four scenarios, resulting in different values for $Q$, with the synchronous dual ascent \eqref{eq:synchronousdai}, in \cref{fig:simulation}.}
The case $Q=1$ corresponds to a synchronous algorithm, where all the agents update their variables at every iteration. For $Q>1$, the agents  perform updates asynchronously, according to their own clocks. \blue{To compare synchronous and asynchronous implementation, we take into account the overall computation burden, i.e., the average number of updates performed per agent}.
For  large values of $Q$, the upper bounds on the step sizes $\gamma_i$'s in \eqref{eq:gamma_bound} decrease, resulting in slower convergence. 
\blue{ However, the bounds can be conservative. In fact, \cref{algo:2} still converges with step sizes set $100$ times larger than their theoretical upper bounds for $Q= 25, 50, 100$ (but not  for $Q=1$).} 
\section{Conclusion}
The distributed dual ascent retains its convergence properties even if the updates are carried out  completely asynchronously and using delayed  information, provided that  small enough uncoordinated step sizes are chosen. Convergence rates for  primal-dual methods in this general asynchronous scenario are currently unknown. 

\bibliographystyle{IEEEtran}
\bibliography{library}

\end{document}